\documentclass[english,a4paper,11pt]{article}
\usepackage{amsmath,amssymb,graphicx,url}
\newtheorem{Theorem}{Theorem}
\newtheorem{Proposition}[Theorem]{Proposition}

\begin{document}
\title{A single fractal pinwheel tile} 
\author{Christoph Bandt, Dmitry Mekhontsev and Andrei Tetenov}
\maketitle

\begin{abstract}
The pinwheel triangle of Conway and Radin is a standard example for tilings with self-similarity and statistical circular symmetry. Many modifications were constructed, all based on partitions of triangles or rectangles. The fractal example of Frank and Whittaker requires 13 different types of tiles. We present an example of a single tile with fractal boundary and very simple geometric structure which has the same symmetry and spectral properties as the pinwheel triangle.
\end{abstract}

\section{Self-similar tilings}
\begin{figure}[h] \label{F1}
\centerline{\includegraphics[width=0.92\textwidth]{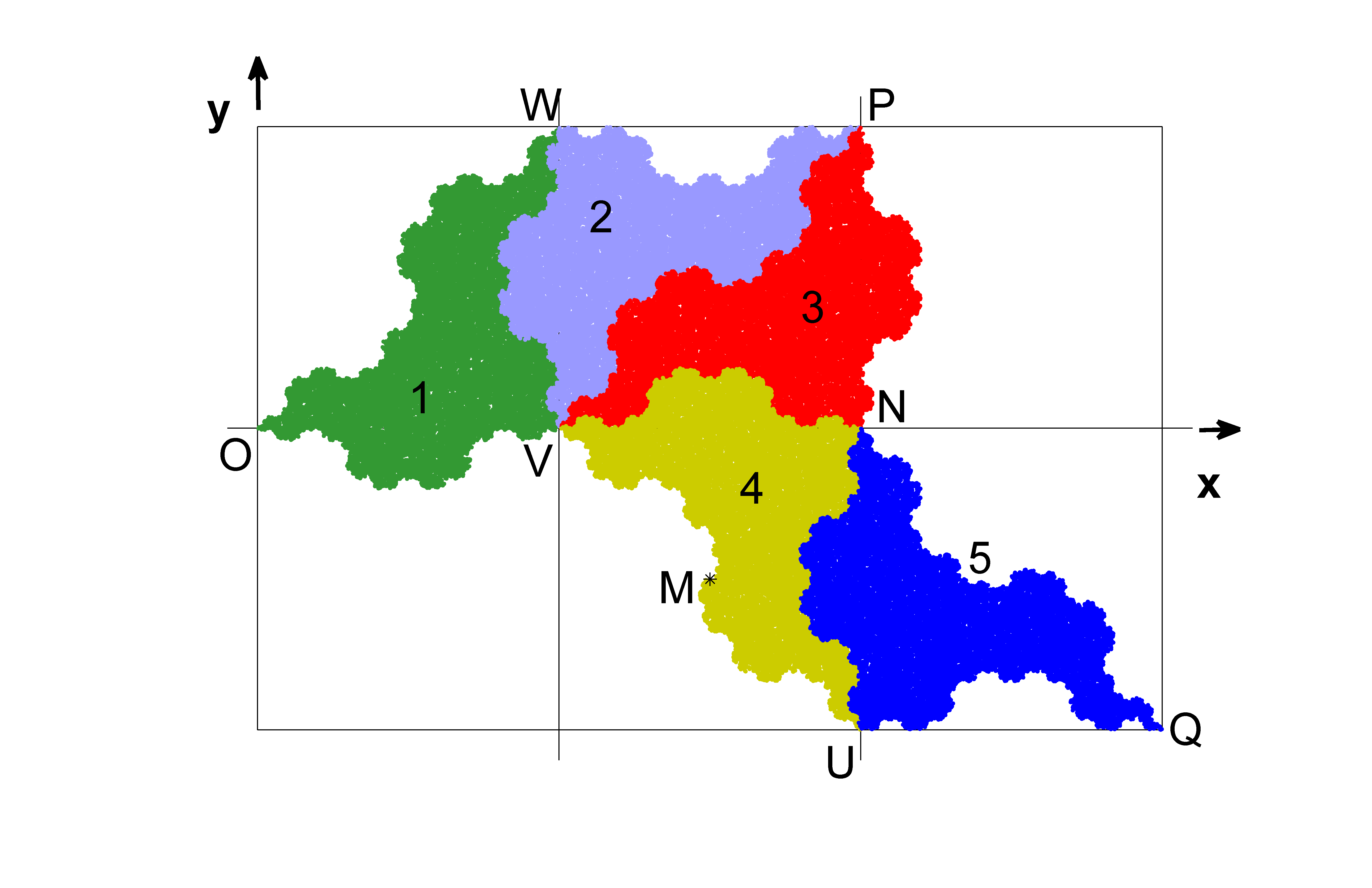}}
\caption{The new fractal pinwheel tile. Numbers $i=1,...,5$ mark the pieces $h_i(A)$ in \eqref{reptiledef}. Capital letters denote landmark points discussed in Section 2.}
\end{figure}

A compact set $A\subset \mathbb R^d$ with non-empty interior is called a \emph{replication tile} with $m$ pieces, or \emph{reptile} for short, if there exist a similarity map $g$ and isometries $h_1,...,h_m$ on Euclidean $\mathbb R^d$ such that
\begin{equation}\label{reptiledef} 
g(A)=h_1(A)\cup ...\cup h_m(A) \, ,
\end{equation}
and any two different sets $h_i(A)$ have no interior points in common. Figure \ref{F1} shows an example, others can be found in  \cite[Chapter 11]{GS}, or in \cite{GG, L,Se,M}.  A well-known theorem of Hutchinson says that $A$ is determined by the data $g,h_1,...,h_m$ \cite{Bar,Fal}.

The union of $m$ isometric copies of $A$ is an enlarged copy $B=g(A),$ so the union of $m$ isometric copies of $B$ will be a still larger copy of $A.$ Continuing this procedure and subdividing the larger copies, one can see that the whole $\mathbb R^d$ is tiled by isometric copies of $A,$ forming a  \emph{self-similar tiling.} There are different ways to form supertiles, by taking the small tile as first, second, or $m$th piece of the larger tile, and one has to care a bit so that this so-called inflation process takes place around an interior point of $A.$ 

Our assumptions are quite restrictive. Among others, comparison of the volume in equation \eqref{reptiledef} shows that all eigenvalues of $g$ must have modulus $\sqrt[d]{m}.$ There are more general concepts in the literature, allowing for affine maps $g,$ for copies $h_i(A)$ with different sizes, and for several types of tiles. Here we stick to the simplest case, and we take $d=2$ so that $g$ has ratio  $\sqrt{m}.$ Moreover, we assume that $h_1=id,$ and consider only tiles $A$ which are homeomorphic to a disk, bounded by a closed Jordan curve. Our point is that even under such restrictive assumptions new examples can be found.

To see why Figure \ref{F1} is important, we give a brief review of self-similar tilings. The main point is that the tiling property of $A$ is fulfilled only for very special choices of similarities $g$ and isometries $h_i.$ The standard assumption is that the mappings $h_i$ generate a certain crystallographic group $\Gamma ,$ and the expanding map $g$ produces a subgroup $\Gamma_0 = g\Gamma g^{-1}\subset \Gamma .$ 

There are tilings which have $\Gamma$ as their symmetry group. This happens when  $H=\{ h_1,...,h_m\}$ is a complete residue system, that is $H\cdot\Gamma_0=\Gamma $ \cite{GG}, cf. \cite{Ba5,Ba97,L}. This necessary and sufficient condition for crystallographic tilings is easy to check. For given data $g,h_i$ there exists only one tiling, and many ways to assemble the tiles into larger and larger supertiles. Examples include the plane regular tilings by squares and by equilateral triangles, with $m=4.$

A second class of tilings with $\Gamma_0 \subseteq\Gamma$ and non-complete residue system has been considered by many authors \cite{LW,GS,Se,Ba97,So,L,Baake2013,LLR}. We have more tilings, but less symmetric tilings, and fewer choices for forming supertiles than in the crystallographic case. The simplest case is $g(x)=4x$ on $\mathbb R$ with $h_i(x)=x+v_i,$ where
$v_i=0,1,8,9$ and  $A=[0,1]\cup [2,3]$ \cite{LW}. There are two possible tilings with a fixed tile $A$ obtained by considering $A$ as left or as right part in a supertile. 
With the IFSTile program package \cite{M}, a lot of new cases in this class have been detected.  Necessary and sufficient conditions for the second class were found only in the one-dimensional case \cite{LLR}.

An extreme subclass of this class includes the chair and sphinx tilings \cite{GS,Se,Baake2013}.  These tilings are not periodic: no translation will transform them into themselves. We have a continuum of different tilings, which can be made a tiling space \cite{SoD,SaB}, and the composition of supertiles is unique in each given tiling \cite{So}. Such tilings are often associated with quasicrystals, although their Fourier spectra are not much different from those of periodic tilings \cite{Se,Baake2013}. The 'typical quasicrystal' Penrose, Robinson and Ammann tilings made from several types of tiles \cite[Chapter 11]{GS} and are not considered here. They have symmetries of the Fourier spectrum which are forbidden for crystallographic pattern, and which have been found by physicists in quasicrystallic alloys \cite{Se,Baake2013}.  

There is a third class of self-similar tilings, where either $H$ does not generate a crystallographic group, or $\Gamma_0 = g\Gamma g^{-1}$ is not a subset of $\Gamma .$ In the language of Section 3 below, the neighbor maps are not contained in a crystallographic group. The first example was the pinwheel triangle of Conway and Radin where it is easy to verify that $\Gamma_0$ contains an irrational rotation. Using ergodicity, Radin \cite{R0} proved that the orientations of triangles in an infinite pinwheel triangle tiling are equidistributed on the circle. This implies that the tiling has a continuous spectrum. A physical material modelled by pinwheel triangles would have an extraordinary diffraction pattern consisting of circles, like a disordered system, cf. \cite{BaaFG,MPS}. At the same time, the pinwheel triangle tilings have a finite set of matching rules, similar to a crystal, as proved in \cite{R}. This apparent contradiction motivated the work of the mathematical physicist Radin. For mathematical work on tiling spaces and spectra of tilings see \cite{SoD,SaB}.

Various modifications of the pinwheel triangle have been presented \cite{Sa, F}. They use triangles, and most of them have a larger number of tiles. There is also a fractal pinwheel version \cite{FW} which uses 13 different types of tiles. 

Figure \ref{F1} shows an unexpected \emph{single fractal tile with irrational rotations between neighbor tiles,} denoted as 'fractal pinwheel'.  Irrational rotations imply statistical circular symmetry of orientations \cite[Section 6]{F}, and continuous diffraction and dynamical spectrum \cite{BaaFG,MPS,SoD}. The assumption of Goodman-Strauss \cite{Go} are also fulfilled, so  there exists a finite set of matching rules. Straight line boundaries are not necessary for such tiles. 

\begin{figure}[h] 
\includegraphics[width=0.999\textwidth]{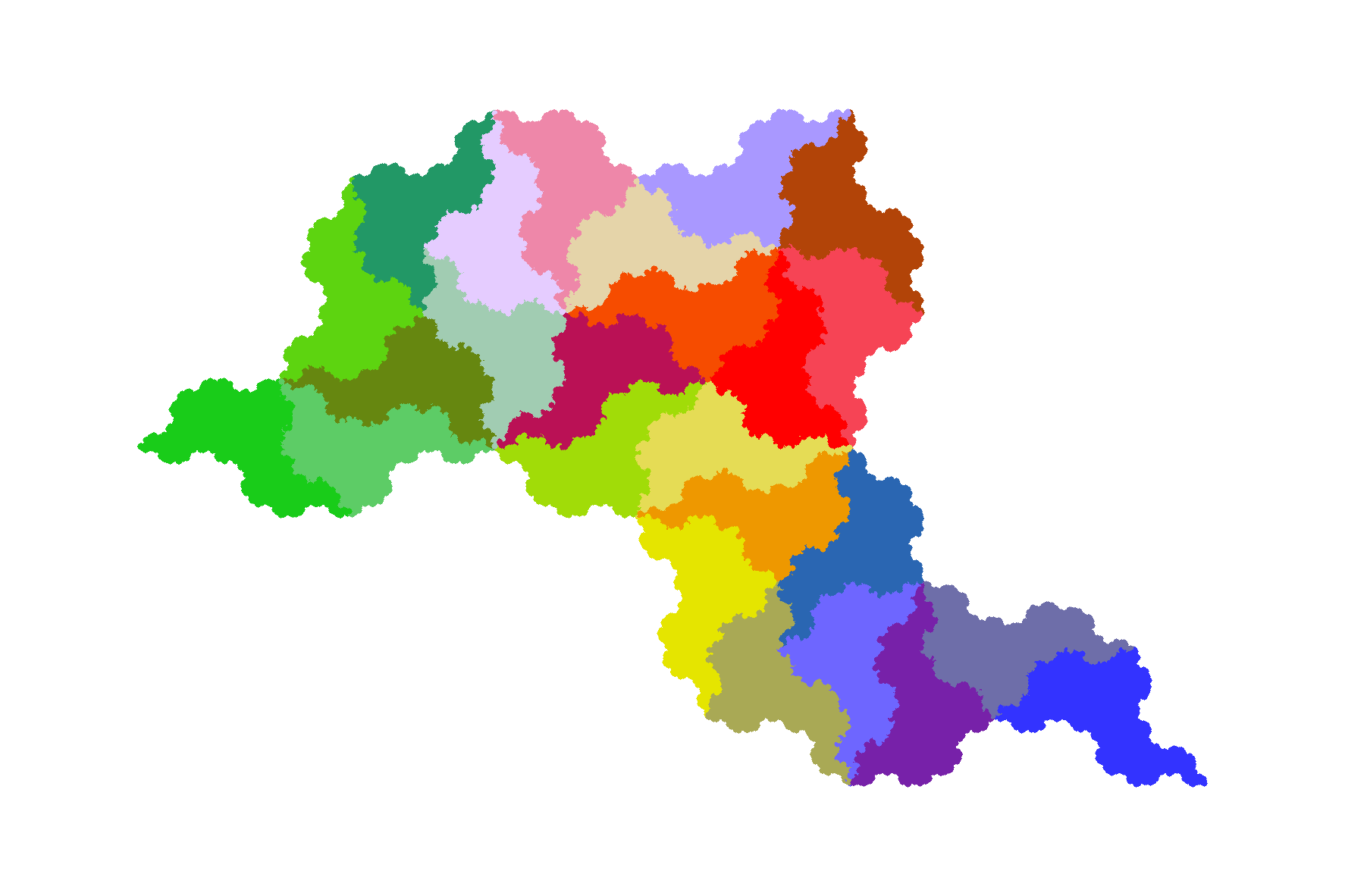}
\caption{Second subdivision of the fractal pinwheel, based on colors as in Figure \ref{F1}. The small piece $h_i(f_j(A))$ is colored with a mixture of 2/3 color $i$ and 1/3 color $j.$ Neighboring subpieces at the border between pieces 1,2 (and 4,5) differ by irrational rotations.} \label{F2}
\end{figure}

In the following section the reptile in Figure \ref{F1} is defined, and geometrical properties are studied. While the tiling property is obvious for a triangle like Radin's pinwheel, considerable effort is needed for the proof in the fractal case where the tile is defined only implicitly by contraction maps. Theorem \ref{T} in Section 3 uses the technique of neighbor maps.  
Actually, the tile was found by a computer search with IFSTile \cite{M} which analyzed neighbor graphs of many random parameter sets. Since then, a search of $10^9$ parameter sets with IFSTile \cite{M} has not provided further pinwheels. 
In Section 4 we study a second reptile structure on the fractal pinwheel $A$ which is quite different.  

After this paper was completed, we found that both of our tiles were already 2012 presented by Ventrella in his inspiring book
\cite[p.85-86]{Ventrella2012} on plane-filling curves. Ventrella used L-systems and gave no rigorous argument for the plane-filling property. He noted that these two curves do not fit into a square grid and show an extraordinary 'mixture of 90 and 45-degree angles'. Our careful analysis of irrational angles will provide mathematical clarity.

\section{Definition and geometry of the fractal tile}
The two outer pieces 1 and 5 of Figure \ref{F1} have different orientation from the whole figure while the three middle pieces have the same orientation. The apparent vertices of the tiles are on a lattice. We take part 1 as our basic tile $A$ and choose coordinates so that $A$ has vertices $O={0\choose 0}, V={1\choose 0},$ and $W={1\choose 1}.$ Then $g(A)$ is a triangle with fractal boundaries and vertices $O, P={2\choose 1},$ and $Q={3\choose -1}.$  For an affine mapping $h{x\choose y}={ax+by+c\choose dx+ey+f},$ the coefficients are given as ${c\choose f}=h{0\choose 0}, \, {a\choose d}=h{1\choose 0}-h{0\choose 0},$ and ${b\choose e}=h{1\choose 1}-h{1\choose 0}.$  
Now $g$ is the linear mapping determined by $g(V)=P, g(W)=Q.$ The $h_i$ are symmetries of the square lattice determined by similar equations, for instance $h_2(O)=V, h_2(V)=W, h_2(W)=P.$ Calculation yields
\begin{eqnarray}\label{coord}   
g{x\choose y}={2 \ \, 1 \choose 1\,  -\! 2}{x\choose y}={2x+y\choose x-2y} \ ,\hspace{15ex} \\   
h_2{x\choose y}={y+1\choose x} \, ,\  
h_3{x\choose y}={2-y\choose 1-x} \, ,\  h_4{x\choose y}={x+1\choose -y}
\, ,\  h_5{x\choose y}={y+2\choose -x}\, .   \label{coor}
\end{eqnarray}

Of course, other descriptions are possible. Choosing part 4 as basic tile $A,$ and point $M$ as origin of a coordinate system with  $N={1\choose 0}, V={0\choose 1},$  we would get the expanding matrix as ${2 -\! 1 \choose  1\  2}.$  Thus the expansion map is the same as for the pinwheel triangle in \cite{R0,R,FW}. The $h_i$ are quite different: in Figure \ref{F1}, piece $k$ is connected only with pieces $k-1, k+1$ by a fractal edge, while in the pinwheel triangle three pieces are pairwise connected by edges or  'half-edges'.  

\begin{Proposition}[Convex hull and diameter of fractal pinwheel] \label{P0}\hfill\\
The convex hull of the fractal $g(A)$ in Figure \ref{F1} is the hexagon $OUQPWL$ with $L={3/5\choose 4/5}\, .$
The diameter of $g(A)$ is $\sqrt{10},$ and the diameter of $A$ is $\sqrt{2},$
\end{Proposition}

{\it Proof. } We denote the convex hexagon $OUQPWL$ by $C.$
The self-similar set $A$ contains the points $O$ and $W$ which are fixed points of mappings $f_1$ and $f_5,$ see \eqref{selfsim} below. Thus $g(A)$ contains $O,W, \, P=h_3(O), Q=h_5(W), U=h_4(W),$ and $L=g^-1(U).$
So the hexagon $C$ is a subset of the convex hull of $g(A).$ On the other hand, it is easy to check that $C_1=g^{-1}(C)\subset H$ and $C_k=h_k(C_1)\subset C.$ This implies that $g(A)\subset C$ (the Hutchinson operator $F$ for $g(A)$ fulfils $F(C)\subset C,$ its iteration yields a decreasing sequence converging to $g(A)$, see \cite{Bar}). Thus $C$ is the convex hull of $g(A).$ It has three sides of length 1, two of length $\sqrt{5},$ and a short side of length $1/\sqrt{5}.$ 

Since the diameter of a polygon is its longest side or diagonal, the diameter of $C$ and of $g(A)$ is $\sqrt{10},$ the length of $OQ.$ Since $g$ has factor $\sqrt{5},$ the diameter of $A$ is $\sqrt{2}.$\hfill $\Box$ \medskip

Two Jordan arcs $\gamma$ and $\gamma '$ with a common endpoint $P$ will be said to form a Jordan angle of size $\alpha$ if a rotation by $\alpha$ with center $P$ transforms $\gamma$ into $\gamma '.$ The fractal boundary curves between $O,P,$ and $Q$ are called sides of the triangle $\triangle OPQ.$

\begin{Proposition}[The fractal pinwheel as a triangle] \label{P1}\hfill
\begin{enumerate}
\item[a)] The fractal triangle $g(A)=\triangle OPQ$ has two congruent sides $OP, PQ.$ The fractal curves $OU$ and $VQ$ on the long side $OQ$ are also congruent to the short sides.
\item[b)] The short sides are symmetric with reflection at the perpendicular bisectors of $OP, PQ.$ The long side is symmetric with respect to a $180^o$ rotation around its midpoint $M.$
\item[c)] The triangle has Jordan angle $90^o$ at $P,$ and irrational angles  $\alpha =\angle UOP=2\cdot\arctan\frac12\approx 53.1^o$ at $O$ and $90^o-\alpha\approx 36.9^o$ at $Q.$
\item[d)] In dimension $\delta =\log (1+\sqrt{2})/\log \sqrt{5}.$ all sides have positive and finite Hausdorff measure  Taking this as length measure, Pythagoras' theorem holds for the triangle.
\item[e)] The area of $A=\triangle OVW$ is $\frac12 .$  
\end{enumerate}
\end{Proposition}

{\it First part of proof. }
Properties a) and b) seem obvious from Figure \ref{F1}, and a proof is given in Section 3. c) immediately follows from a). The Euclidean triangle $\triangle OPN$ is a Radin pinwheel and has angle $\arctan\frac12$ at $O.$ The size of the Jordan angle $\angle VQP$ can be checked by noting that $\angle VQZ=90^o,$ where $Z={4\choose 1}$ is outside Figure \ref{F1}. Irrational angles are the basis for concluding that the tilings are statistically circular symmetric \cite{R,F}.

Concerning d), we shall prove in Section 3 that the sides of $\triangle OPQ$ form a graph-directed system of self-similar sets. Calculation of the Hausdorff dimension $\delta$ of such fractal boundaries is standard, starting with classical work of Gilbert in the 1980s \cite{Gi,DKV,SW,HLR,BM}. The Hausdorff measure $\mu$ of dimension $\delta$ of each side is positive and finite. Assuming this fact, we can prove d) by the following argument.

Hausdorff measure is invariant under isometry, and Hausdorff measure of an image under a similitude with factor $r<1$ decreases by the factor $r^\delta .$ Here the factor of $g^{-1}$ is $r=1/\sqrt{5},$ and we let $y=r^\delta .$ If we put $a=\mu (OP)=\mu (PQ)$ then
\[  \mu(OQ)=\mu(OU)+\mu(UQ) =a(1+y)\ \mbox{ and }\  a=\mu(OW)+\mu(WP)=y(\mu(OQ)+a) \ . \]
Now $a$ cancels out, and $y=\sqrt{2}-1$ is the positive solution of the quadratic equation. Thus $ \mu(OQ)=a\sqrt{2}$ which means that Pythagoras' theorem is true for our triangle with fractal side lengths. Under contractions like $g^{-1}$ where Euclidean distances shrink by $r=1/\sqrt{5}\approx .447 ,$ the lengths of fractal boundaries shrink faster, by $\sqrt{2}-1\approx .414.$ The dimension of boundaries is $\delta=\log y/\log r \approx 1.1\, ,$ as stated in d).

It is a well-known open problem whether there exists a single puzzle tile $T$ such that the plane can be tiled with isometric copies of $T,$ but it cannot be tiled periodically \cite[Chapter 11]{GS}.  Could Figure \ref{F1} be such an aperiodic tile? Unfortunately not. If we add to $g(A)$ three copies obtained by successive $90^o$ rotation around $P,$ we get a fractal square which tiles the plane like an ordinary square. This directly follows from a) and b).  Another periodic tile is the fractal rectangle $PNVW,$ the union of parts 2 and 3 in Figure \ref{F1}.  If we consider the square lattice as a checkerboard, and we put a copy of the rectangle on each white square, and a copy rotated around $90^o$ on each black square, we have a tiling, due to Proposition \ref{P1} a) and b). Part of such checkerboard pattern can be seen in Figure \ref{F2}.
As a consequence, the fractal rectangle must have area one, and the area of the tile $A$ is $\frac12 .$ This proves e). 
\hfill $\Box$ \medskip

\begin{figure}[h]
\includegraphics[width=0.999\textwidth]{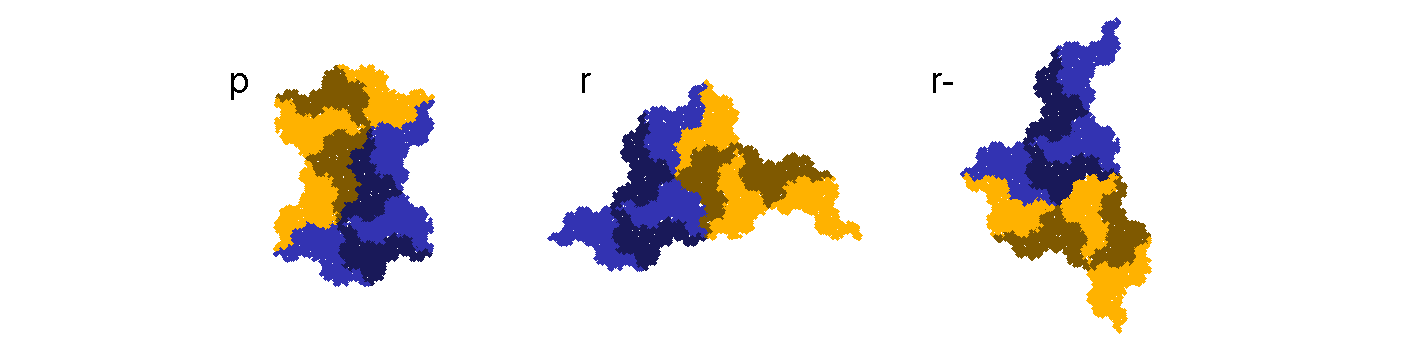}\\
\includegraphics[width=0.999\textwidth]{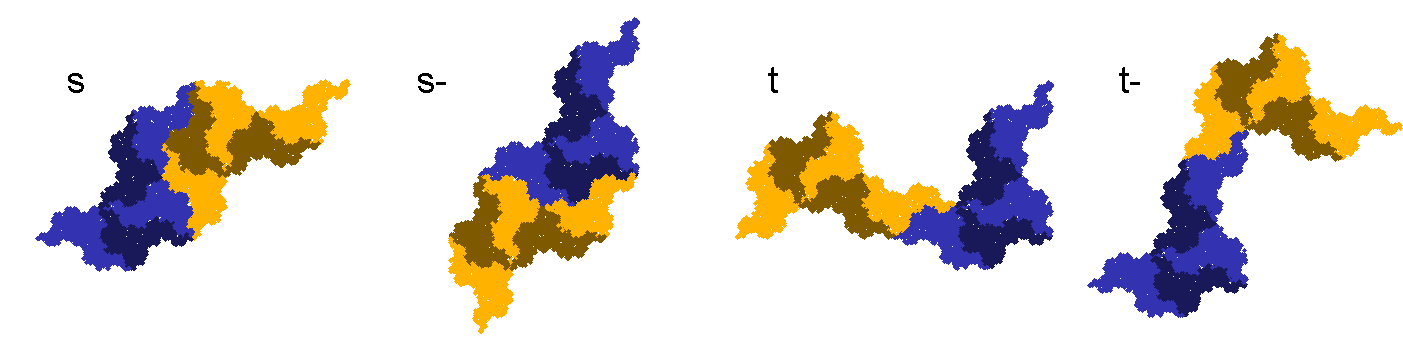}\\
\includegraphics[width=0.999\textwidth]{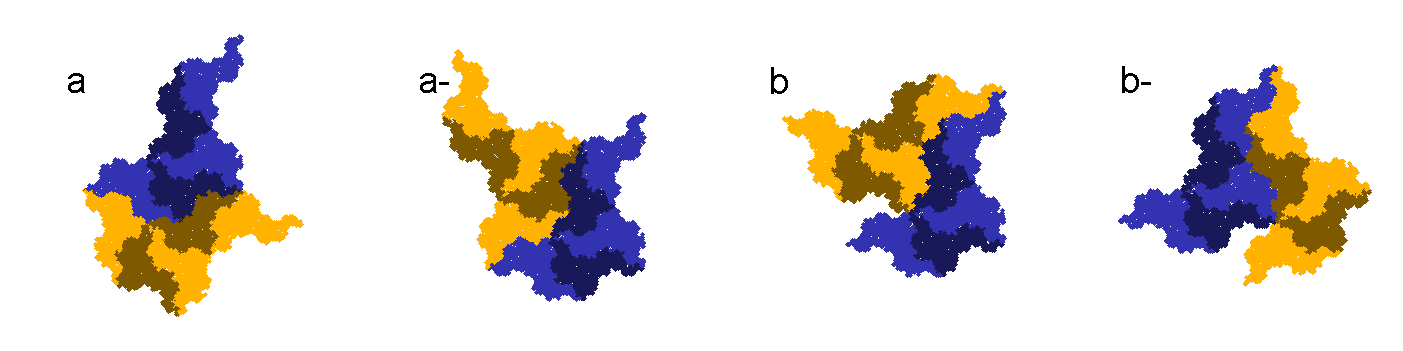}
\caption{Edge neighbors. Pieces number 2 and 4 are displayed darker in $A$ and in the neighbor $h(A)$ to recognize neighbor pieces which define arrows in the neighbor graph.}\label{F3}
\end{figure} 

\section{The neighbor graph}
Applying $g^{-1}$ to both sides of \eqref{reptiledef}, we represent $A=A_1$ as a self-similar set or IFS attractor
 \begin{equation}\label{selfsim} 
A=f_1(A)\cup ...\cup f_m(A) \quad\mbox{ with }\quad f_i=g^{-1}\cdot h_i   \ \mbox{ for } i=1,...,m.
\end{equation} 
For our example, let us list the expressions of $f_i{x\choose y}, i=1,..,5$ to be used below.
\begin{equation}\label{fi} 
\frac15{ 2x+y\choose x-2y},\ \frac15{x+2y+2\choose -\! 2x+y+1}, \ \frac15{-\! x -\! 2y+5\choose\! 2x - y\ }, \
\frac15{2x -\! y+2\choose x+2y+1}, \ \frac15{-\! x+2y+4\choose 2x+y+2} .
\end{equation}
A basic theorem of Hutchinson says that there is exactly one compact solution $A$ of \eqref{selfsim} for given contractions $f_i$  \cite{Bar,Fal}. To get good topological and geometric properties, however, we have to assume that the $f_i(A)$ do not overlap too much. Usually one assumes the open set condition: there is an open set $U$ so that the $f_i(U)$ are disjoint subsets of $U.$ Such a set $U$ is hard to determine, so one requires that $A$, or better $f_1,...,f_m,$ is of finite type. This means that the neighbor graph, defined below, is finite.  The topology of $A$ is determined by a finite automaton which is obtained from the $f_i$ by an algorithm. We shall calculate the automaton for our example.

The technique of neighbor graphs is now well established. The beginnings go back to Gilbert in the eighties, see \cite{Gi} and references in \cite{Ba5}. In \cite{BG}, neighbor maps were introduced as $h=f_{\bf i}^{-1} f_{\bf j}$ where ${\bf i,j}\in \{ 1,...,m\}^n$ for some $n\ge 1.$ The open set condition was shown  to be equivalent to the fact that neighbor maps cannot converge to the identity map. For tiles $A,$ neighbor maps are exactly the isometries between neighboring pieces in any self-similar tiling made up of copies of $A,$ cf. \cite{Ba97}. We shall consider only proper neighbors, for which  $f_{\bf i}(A)\cap f_{\bf j}(A)\not=\emptyset .$  Here $\bf i= i_1...i_n$ and $f_{\bf i}= f_{i_1}\cdot...\cdot f_{i_n}.$ Thus a neighbor map transforms $A$ to an isometric copy $h(A)$ which intersects $A$ and is a tile in some patch obtained by inflation of $A$ around one of its pieces $f_{\bf i}(A).$

A neighbor map $h=f_{\bf i}^{-1} f_{\bf j}$ describes a potential boundary set $A\cap h(A).$ Actually, the study of tile boundaries, especially for the Levy curve \cite{DKV,SW} was a key motivation for developing the method of neighbor graphs \cite{HLR,ScT,BM,AL}. The neighbor graph yields a system of set equations for the boundary of $A,$ similar to \eqref{selfsim}. In the program IFSTile \cite{M}, a fast and very general algorithm was implemented to determine neighbor graphs and dimensions of boundary tiles.

\begin{figure}[h]  \hspace*{7ex}
\includegraphics[width=0.4\textwidth]{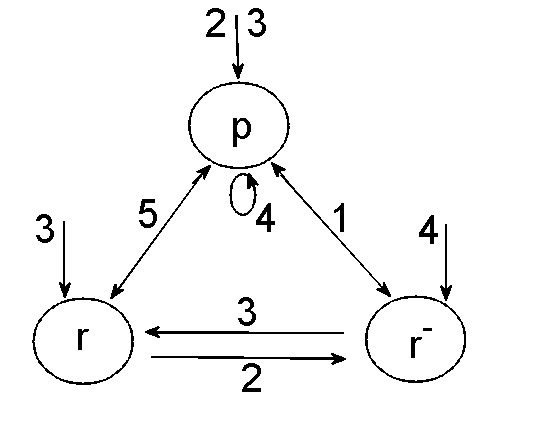}\hspace{7ex}
\includegraphics[width=0.4\textwidth]{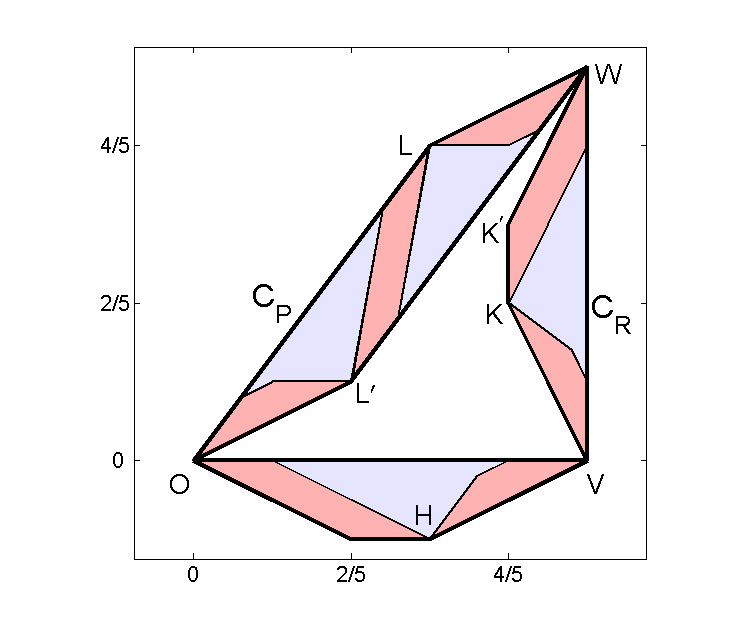}
\caption{Left: graph of edged neighbors by rational rotation, arrows marked by first label only. Right: Convex hulls of associated boundary sets and their images in \eqref{quadryeq}.} \label{G1}
\end{figure}

The vertices of the neighbor graph $G$ are the neighbor maps $h.$ An arrow with label $i,j$  is drawn from vertex $h$ to vertex $\bar{h}$ if $\bar{h}=f_i^{-1}hf_j,$ for two marks $i,j\in\{ 1,...,m\} .$  We keep only those arrows which correspond to proper neighbors, that is, $f_i(A)\cap h(f_j(A))\not=\emptyset.$ See \cite{BM} for details.  The identity map $id$ is the root vertex of the graph, with loops labelled $i,i.$ It is not drawn in Figure \ref{G1}, and arrows from $id$ have no intial vertex. 

If $G$ is a finite graph, the reptile or IFS attractor $A$ generated by $g,h_1,...,h_m$ resp. $f_1,...,f_m$ is called \emph{finite type.} If there are no incoming edges to the root vertex $id,$ then the open set condition is fulfilled, which together with the condition that all $f_i$ have similarity ratio $1/\sqrt{m}$ implies the tiling property \cite{BG,Ba97}. 

In a plane tiling, we can have two different kinds of neighbors:  \emph{point neighbors} which have a single intersection point, and  \emph{edge neighbors} which have uncountably many points in common. Other kinds of neighbors can occur  \cite{BG} but not in the case of our example. Moreover, the intersection of two edge neighbors is always homeomorphic to an interval, as will be proved now.

\begin{Theorem}[Neighbor graph and boundary of fractal pinwheel] \label{2D}\hfill\\ \label{T}
Let $A$ denote the fractal pinwheel, with mappings defined by \eqref{selfsim}, \eqref{coord} and \eqref{coor}.
\begin{enumerate}
\item[a)] There are exactly 11 edge neighbors illustrated in Figure \ref{F3}, 69 point neighbors and no other neighbors. Thus $A$ is finite type, has non-empty interior and is a reptile.
\item[b)] Two of the maps for edge neighbors are irrational rotations. So there is a continuum of different tilings. They are not lattice tilings, and have statistical circular symmetry.
\item[c)]  Edges are of two types: $k\cdot 90^o$ rotation on one hand, glide reflections and irrational rotations on the other. All subedges of an edge at any level have the same type as the original edge.
\item[d)] $A$ is homeomorphic to a disk, bounded by a closed Jordan curve of dimension $\delta =\log (1+\sqrt{2})/\log \sqrt{5}.$ This boundary set is the union of intersections of $A$ with its three neighbors by rational rotation.
\end{enumerate}
\end{Theorem}

{\it Proof. }First we sketch our proof of the difficult part a) by calculation of all neighbor maps with computer. This was done independently by two authors with different software.  We build the graph $G$ recursively, calculating all possible $h'=f_i^{-1}hf_j$ for all previously constructed maps $h.$   We want to neglect $h$ if $A\cap h(A)=\emptyset$ but this cannot be checked directly.  However, by Proposition \ref{P0}, the diameter of $A\cup h(A)$ is smaller than $2\sqrt{2}$ whenever $A$ intersects $h(A).$  Since the origin belongs to $A,$ this implies that $\| h{0\choose 0}\|\le\sqrt{8}.$

We determine the graph of neighbor maps $h{x\choose y}=O{x\choose y}+{e\choose f}$ with an orthogonal matrix $O$ for which $e^2+f^2\le 8$ is fulfilled.  
This graph turns out to be finite, with 955 vertices. Then we take the subgraph of all vertices which lie on cycles of the large graph. This is our graph $G$ of proper neighbors with only 81 vertices including the root. Point neighbor maps $h$ have the property that for each $n,$ only one path of length $n$ starts at vertex $h.$ They are easily singled out by checking powers of the adjacency matrix of $G.$ There were 69 point neighbors, 11 remaining neighbor maps and the identity, which proves a). On a PC, all this is done in less than 2 seconds.

Now we give a computer-free proof of the theorem, except for the number of point neighbors. As explained below, the 11 edge neighbors in Figure \ref{F3} can be found by inspection of the second subdivision of $g(A),$ Figure \ref{F2}, and confirming by calculation. A check of the next subdivision, or of Figure \ref{F3}, then verifies that no other edge neighbors exist, and there are no incoming arrows to the root vertex in the neighbor graph.  

Note that \eqref{selfsim} implies $f_i^{-1}=h_i^{-1}g$ and thus $f_i^{-1}f_j=h_i^{-1}h_j$ for $i,j=1,...,m=5.$ So the successors of the root vertex $id$ include the rational rotations $p=h_2^{-1}h_3=h_3^{-1}h_2=-x+{1\choose 1},$ a $180^o$ rotation around ${\frac12\choose\frac12}$ and $r=h_3^{-1}h_4,$ a clockwise $90^o$ rotation with center 
$V={1\choose 0}.$  Since the inverse of any neighbor map $h=f_{\bf i}^{-1} f_{\bf j}$ is the neighbor map 
$f_{\bf j}^{-1} f_{\bf i},$ we have to add the inverse $r^{-}=h_4^{-1}h_3,$ the counterclockwise $90^o$ rotation with center $V.$ The mapping $p$ is self-inverse.

So far we have studied the maps between pieces 2 and  3, and 3 and 4 in Figure \ref{F1}. Now we consider their subpiece neighbors in Figures \ref{F2} or \ref{F3}.  We see that subpieces 24 and 34 have the same relative position as pieces 2 and 3, which is algebraically verified by the equation $p=f_4^{-1}pf_4$ and results in a loop from vertex $p$ to itself with label $4,4.$
Subpieces 41 and 35 also intersect, and correspond to the neighbor map $p,$ which results in arrows from vertex $r$ to $p$ with label $1,5,$ and from $r^{-}$ to $p$ with label $5,1.$ Checking two other pairs of subpieces of 2 and 3, and one remaining pair of subpieces of 3 and 4 in Figure \ref{F2}, we obtain the graph in Figure \ref{G1}. This argument proves that there are no other arrows starting in $p,r,r^-$ (which the computer checked algebraically).  It is enough to take only the first label $i$ of any arrow from a vertex $h$ to a vertex $\bar{h}$ since the second label $j$ is the same as the first label of the arrow from $h^{-1}$ to $\bar{h}^{-1}.$ Drawing arrows from the root without an initial vertex, we obtain a reduced form of the graph \cite{BM} on the left of Figure  \ref{G1}.

\begin{table}\label{arrows}
\tabcolsep1.5mm
\quad\begin{tabular}{|l|cc|c|ccc||cc|ccc|ccc|c|ccc|}\hline
&\multicolumn{6}{c||}{rational rotations}& \multicolumn{12}{c|}{glide reflections, irrational rotations}\\ \hline
initial vertex&\it id&\it id&\it p&\it p&\it r&\it r&\it id&\it id&\it s&\it s&\it s&\it a&\it a&\it a&\it t&\it b&\it b&\it b\\ \hline
terminal vertex&\it p&\it r&\it p&\it r&$p$&$r^{-}$&\it s&\it s&\it a&\it t&\it b&$a^{-}$&$t^{-}$&$b^{-}$&$s^{-}$&$a^{-}$&$t^{-}$&$b^{-}$\\ \hline
first label&2&3&4&5&5&3&1&4&3&5&5&1&1&2&1&4&4&5\\ \hline
second label&3&4&4&1&1&5&2&5&1&1&2&1&4&4&5&3&5&5\\ \hline 
\end{tabular} \vspace{2ex} 
\caption{Arrows between edge neighbors and their labels. To each arrow $(h,h',i,j)$ there is another arrow 
$(h^{-1},h'^{-1},j,i)$ which is not listed here, except for $(p,p,4,4)$ with $p=p^{-1}.$ }
\end{table}

To get all edge neighbors, we still have to consider the boundary between pieces 1 and 2, or 4 and 5.  There we get the glide reflection $s=h_2=h_4^{-1}h_5$ and its inverse $s^{-}$ seen in the second row of Figure \ref{F3}. The subpieces 15 and 21 lead to the glide reflection $t=f_{15}^{-1}f_{21}=f_5^{-1}s$ and its inverse $t^{-}.$ Subpieces 13 and 21 yield the neighbor map $a=f_{13}^{-1}f_{21}=f_3^{-1}s$ which is an irrational rotation by $\alpha$ around $O,$ see Section 2 and Figure \ref{F3}.  Subpieces 15 and 22 yield the neighbor map $b=f_{15}^{-1}f_{22}=f_5^{-1}sgf_2$ which is an irrational rotation by $90^o-\alpha$ around $W={1\choose 1}.$ 

We found edge neighbors in the second subdivision for which the neighbor map is an irrational rotation! This shows the non-crystallographic character of our fractal tile.  This property implies that there is a continuum of different tilings and that for each tiling, the orientations of tiles, defined as angles, are dense in $[0,2\pi].$  Their distribution within a large circle of radius $R$ around 0 converges to the uniform distribution on $[0,2\pi]$ when $R$ runs to infinity. This is called 'statistical circular symmetry' \cite{F}. The Fourier spectrum, important from the physicists viewpoint, is also symmetric under rotations. This was shown in \cite{R,MPS,F} which completes the proof of  b).

To get the complete graph of edge neighbors, we still have to study the subpieces of neighbors $t,a,b,t^{-},a^{-},b^{-}$ in Figure \ref{F3}. They all represent neighbor maps of the second and third row of  Figure \ref{F3}, providing arrows in $G$ leading to previous vertices. Instead of drawing this part of $G,$ which is not planar, we list the arrows in Table \ref{arrows}. To each arrow $(h,h',i,j)$ in the table, except $(p,p,4,4),$ there is another arrow 
$(h^{-1},h'^{-1},j,i)$ which is not listed for brevity. This proves a) when we neglect point neighbors. Assertion c) follows since the right part of Table \ref{arrows} contains no arrows leading to rational rotations. The graph of edge neighbors, without root, splits into two components. 

Can we really neglect point neighbors? Yes, we can. The proof of d) will be done only with the graph of rational rotations in Figure \ref{G1} which was derived by simple calculation. d) implies that the three sides of the triangle $A$ studied in Section 2 are really Jordan curves, and thus the angles $\alpha, 90^o-\alpha$ are correctly defined. Together with the list of edge neighbors and the Jordan curve theorem, this implies that any non-edge neighbor can intersect $A$ only in one of the vertices $O,V,W,g^{-1}(U)$ or $g^{-1}(V).$ Moreover, since subtiles meet at such a point with their vertices, only finitely many angles are possible. This shows that beside edge neighbors, only finitely many point neighbors exist, and shows the finite type property of $A.$ (If we are satisfied with the open set condition for $A,$ instead of finite type, the finite number of angles will not be needed.)

Now let us prove d).  Abstract methods as in \cite{AL,BM} are not needed since our case is rather simple.
Consider $A$ with three edge neighbors defined by $p,r,r^{-}.$  The subpiece 33 in the middle of Figure \ref{F2} has this structure. Let $P=A\cap p(A), R=A\cap r(A)$ and $R^{-}=A\cap r^{-}(A)$ denote the corresponding boundary sets of $A.$  The reduced form of the neighbor graph for this configuration in Figure \ref{G1} yields the equation system
\begin{equation}\label{boundaryeq}
P= f_5(R)\cup f_4(P)\cup f_1(R^{-})\, ,\ R= f_3(R^{-})\cup f_5(P)\, ,\ R^{-}=f_1(P)\cup  f_2(R)\ .
\end{equation}
The sets $P,R,R^{-}$ form a so-called graph-directed construction. The crucial point is that the convex hulls $C_P,C_R$ and $C_{R^{-}}$ provide the open set condition for this system. Similar to Proposition \ref{P0}, these are the quadrilaterals
$C_P={\rm conv\,} OLWL'$ with $L={3/5\choose 4/5}$ and $L'={2/5\choose 1/5}\, ,$
$C_R={\rm conv\,} WVKK'$ with $K={4/5\choose 2/5}$ and $K'={4/5\choose 3/5}\, ,$ and 
$C_{R^{-}}={\rm conv\,} OVHH'$ with $H={3/5\choose -1/5}$ and $H'={2/5\choose -1/5}\, .$ The open set condition says that the interiors of the quadrilaterals contain disjoint unions of their images defined in \eqref{boundaryeq}:
\begin{equation}\label{quadryeq}
C_P\supset f_5(C_R)\cup f_4(C_P)\cup f_1(C_{R^{-}})\, ,\ C_R\supset f_3(C_{R^{-}})\cup f_5(C_P)\, ,\ 
C_{R^{-}}\supset f_1(C_P)\cup  f_2(C_R)\ .
\end{equation}
This is verified by simply calculating images of vertices with \eqref{fi}. We get chains where each quadrilateral has one vertex in common with its predecessor and successor. Moreover, each of the points $O,V,W$ belongs to two components of the boundary $J=P\cup R\cup  R^{-}.$ (Details: addresses of points of a boundary set are given by the paths starting in the corresponding vertex of the graph in Figure \ref{G1}. Paths with $555...=\overline{5}$ start in both $p$ and $r.$ So the point with address $\overline{5},$ fixed point $W$ of $f_5,$ belongs to $P\cap R.$ Similarly, fixed point $O$ of $f_1$ belongs to $P\cap R^{-}.$ Since $V=f_2(W)=f_3(O)$ and paths labelled $3\overline{1}, 2\overline{5}$ start in $R,R^{-},$ respectively, $\{ V\} =R\cap  R^{-}.$)

The intersection points of consecutive small quadrilaterals belong to $A$ since they are images of such intersection points on previous levels, for example $L'=f_4(O), L=f_4(W).$
Iterating the graph-directed construction on quadrilaterals we obtain longer chains of smaller quadrilaterals  with vertices in $J.$ This is a classical 'Koch curve' construction. In the limit we have three Jordan arcs $P,R,R^{-}$ which form the closed Jordan curve $J.$ By definition, $J\subset A.$

We show that $A$ has no points in the exterior region of $J.$ The neighbors $p(A), r(A),$ $r^{-}(A)$ contain the closed Jordan arcs $p(J), r(J),$ and $r^{-}(J)$ which have similar neighborhoods of quadrilaterals as $J.$ Comparing slopes of lines, we see that the convex hull of $A,$ given as the outer boundary in Figure \ref{G1}, is within the interior region of the union $J\cup p(J)\cup r(J)\cup r^{-}(J)$ (for the neighbors, use inner sides of quadrilaterals as bound). Thus each point of $A$  exterior to $J$ must belong to one neighbor. So by definition it belongs to $J.$ Thus such exterior points cannot exist. 

We note that the fractal arc $R$ is invariant under the reflection $\sigma$ at the line $y=\frac12 .$ To see this, we check that for each point connecting two quadrilaterals in the approximating chain of $R$ on some level $n$, the reflected point will also be on two quadrilaterals, at least on level $n+1.$ For $K$ and $K'$ this can be seen in Figure \ref{G1}. By induction we prove that \emph{all} vertices of the quadrilaterals within $C_R$ lie on $R.$ They form a reflection symmetric set.

As a consequence, the neighbor map $s=\sigma\cdot r$ describes the same boundary set $R$ as $r.$ On one hand $\rho(R)=R$ implies that $s(A)$ contains $R.$ On the other hand, points in $r(J)\setminus R$ fulfil $x>1,$ and do their images under $\rho,$ so that $s(A)$ cannot contain other points of $A.$

Now $J$ and all boundaries between the pieces $f_i(A)$ form a network of Jordan curves which belong to $A,$ because the edge neighbor maps within $g(A)$ are $p,r,r^{-},$ and $s.$ We can apply the $f_i$ to the union of all these Jordan curves and get a more dense network of Jordan curves bounding the second level pieces and forming a subset of $A.$ The diameter of holes within this network is at most $\sqrt{2}/5.$ Iterating further, the diameter of holes tends to zero. Thus the closed set $A$ contains the whole interior region of $J.$

Once we know that $A$ is homeomorphic to a disk, we immediately have the open set condition and the tiling property. For topological reasons, disk-like neighbors can only meet in a Jordan arc or in a single point. So all remaining neighbors of $A$ are point neighbors. As mentioned above, this implies the finite type property. Since we had a graph-directed system \eqref{boundaryeq} with open set condition \eqref{quadryeq}, the calculation of the dimension $\delta$ in Section 2 is justified. The computer-free proof of the theorem is finished.    \hfill $\Box$ \medskip

{\it Completion of proof of Proposition \ref{P1}. } We know that $A$ is homeomorphic to a disk, and the neighbors $r(A), p(A)$ and $r^{-}(A)$ intersect $A$ in the fractal arcs $R=VW, P=WO,$ and $R^{-}=OV,$ respectively. This immediately implies that the long side $P$ is invariant under $180^o$ rotation $p,$ and that  $R^{-}$ is mapped by $90^o$ rotation $r$ onto $R.$ Reflection-invariance of $R$ was shown above.  The congruence of $R$ with the fractal arc $WL',$ and of $R^{-}$ with $OL$ is given by the neighbor maps $b$ and $a,$ see Figure \ref{F3}. This proves the corresponding statements a)--c)  for $g(A).$ The argument of d) was justified above, and e) is based on the neighbor map $p$ and the reflection invariance of $R.$  Everything is proved.
 \hfill $\Box$ \medskip

\begin{figure}[h]  \hspace*{7ex}
\includegraphics[width=0.85\textwidth]{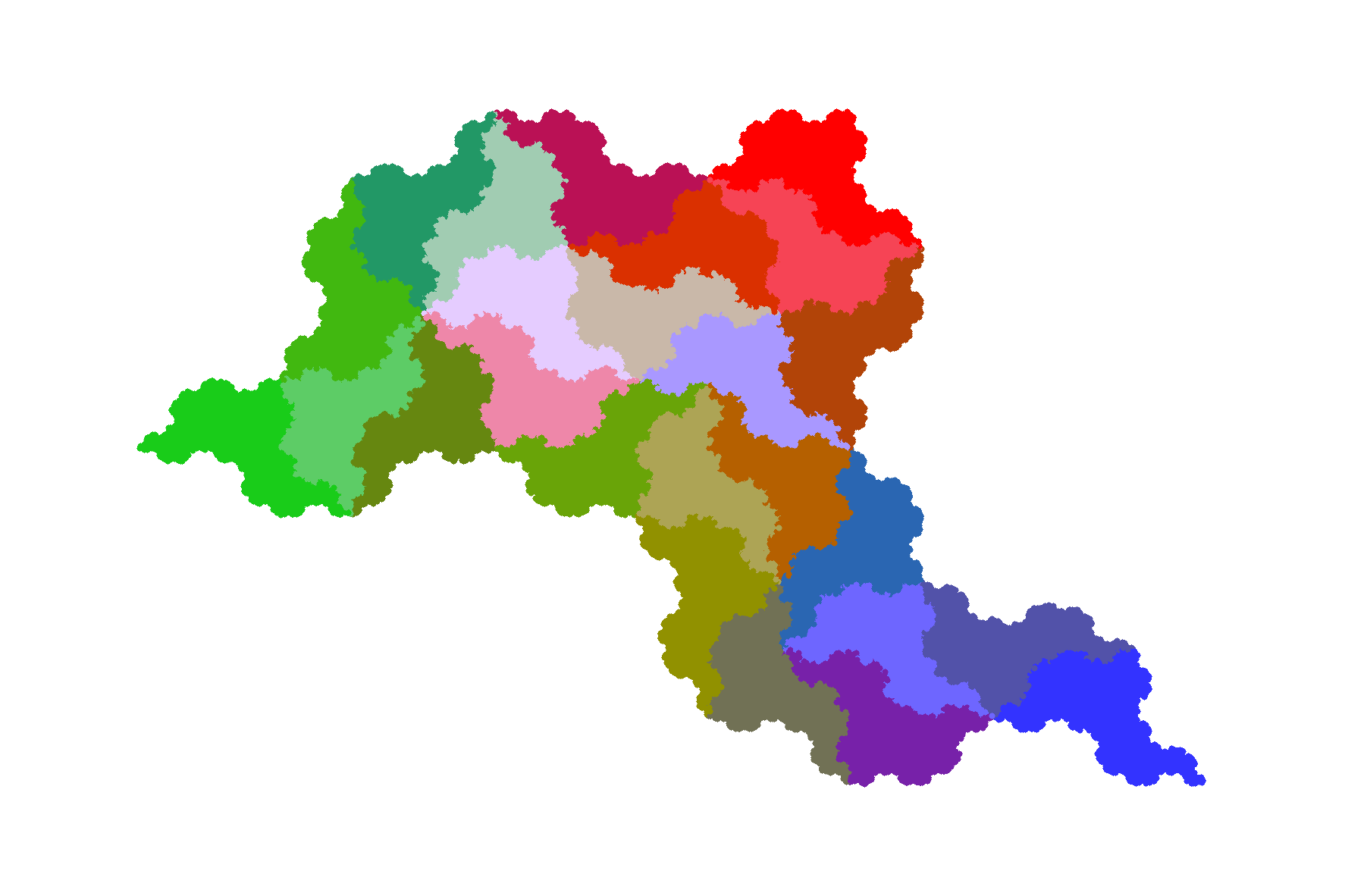}
\caption{Second subdivision of the second fractal pinwheel}\label{F22}
\end{figure}

\section{The second inflation structure}
The union of pieces 2 and 3 in Figure \ref{F1} has the symmetry group of a rectangle. It is mapped to itself by the $180^o$ rotation $p=h_2^{-1}h_3,$ and also by reflection $\sigma {x\choose y}={x\choose 1-\! y}$ at the line $y=\frac12 .$ This kind of symmetry is rare in fractal tiles. We apply the reflection to $h_2$ and $h_3,$ obtaining new maps
\begin{equation}
\bar{h}_2=\sigma h_2={y+1\choose 1- x}\ ,\quad \bar{h}_3=\sigma h_3={2-y\choose  x}\ . \label{h22}
\end{equation}
Since $\bar{h}_2(A)\cup \bar{h}_3(A)=h_2(A)\cup h_3(A),$ there is a new reptile with maps 
$g,h_1,\bar{h}_2, \bar{h}_3, h_4,$ and $h_5.$ As a set, this reptile coincides with $A,$ but the subdivision is different, as shown in Figure \ref{F22}. This leads to other tilings. 

\begin{figure}[h]
\includegraphics[width=0.999\textwidth]{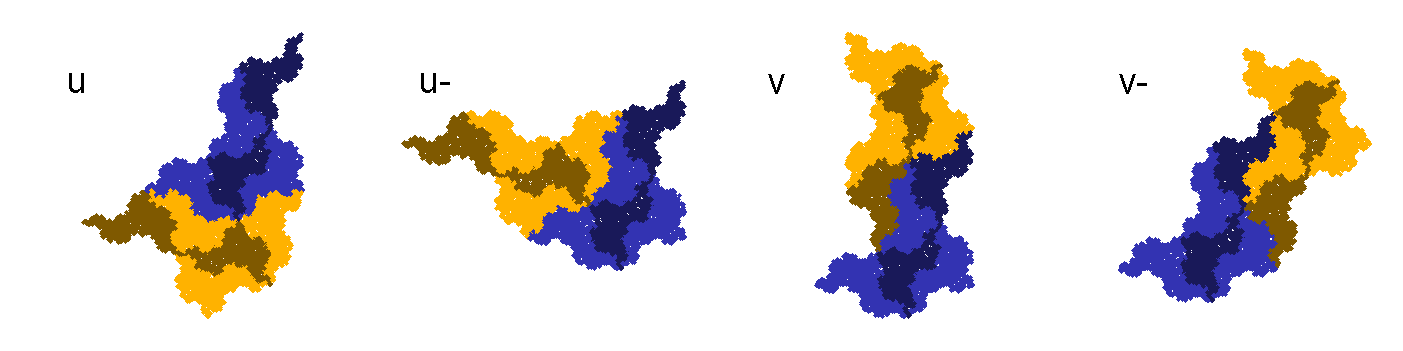}
\caption{For the second structure, these glide reflections replace the irrational rotations in the last row of Figure \ref{F3} as neighbor maps. Pieces 2 and 5 are displayed darker in $A$ and $h(A)$ to recognize adjacent pieces.}
\label{F32}
\end{figure}

Reflection of a rectangle consisting of two pieces in a self-similar triangle was Radin's trick to come from a crystallographic tile to the non-crystallographic pinwheel triangle. A similar trick was used by Conway and Radin \cite{CR} to obtain three-dimensional quaquaversal tilings from crystallographic ones. In our case, however the reflection of pieces 2,3 leads from one non-crystallographic tile to another non-crystallographic tile. 

\begin{figure}[h]
\begin{center}
\includegraphics[width=0.6\textwidth]{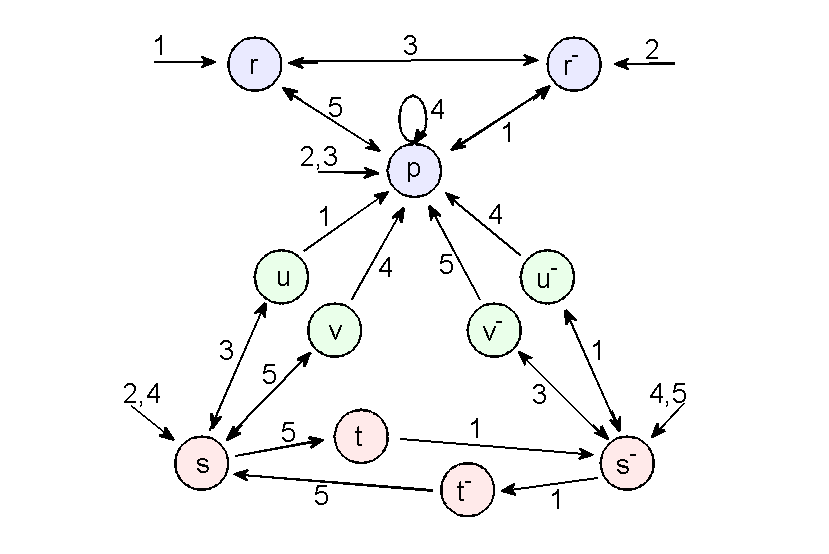}
\caption{Graph of edge neighbors for the second inflation structure}\label{nbg2}
\end{center}
\end{figure}

The second subdivision of this fractal structure is shown in Figure \ref{F22}. Four of the mappings $f_i$ are orientation-reversing. All three vertices of the fractal triangle are fixed points of corresponding contraction maps, resulting in a smaller number of point neighbors. The graph of edge neighbors is planar, as shown in Figure \ref{nbg2}. This second similarity structure has quite different matching rules than the first. The irrational rotations in the last row of Figure \ref{F3} do not appear as neighbors maps. Instead, we have the glide reflections shown in Figure \ref{F32}. An irrational rotation occurs  between point neighbors, as for the pinwheel triangle:  $h=f_{44}^{-1}f_{51}$ has the form 
$h{x\choose y}={.8x-.6y+1\choose .6x+.8y-1}$ with $h(L)=V$ in {\it both} structures. Part a) of the following statement is proved like Theorem \ref{T}. c) follows from the graph in Figure \ref{nbg2}. 

\begin{Proposition}[Neighbor graph of the second fractal pinwheel structure] \label{2D}\hfill\\ \label{T}
The second fractal pinwheel structure on $A$ has the following properties.
\begin{enumerate}
\item[a)] There are exactly 11 edge neighbors, given in the first two rows of Figure \ref{F3} and in Figure \ref{F32}. 
Their graph is shown in Figure \ref{nbg2}. There are 35 point neighbors and no other neighbors. Thus $A$ is finite type, has non-empty interior and is a reptile.
\item[b)] Point neighbor maps include irrational rotations. So there are no lattice tilings, and we have statistical circular symmetry.
\item[c)]  Edges come with two types of maps: rational rotations on one hand, glide reflections on the other. All subedges of an edge of first type are again of first type. Edges of second type contain a dense set of subedges of first type.
\end{enumerate}
\end{Proposition}

\bibliographystyle{plain}
\bibliography{lit0}

\providecommand{\bysame}{\leavevmode\hbox to3em{\hrulefill}\thinspace}
\providecommand{\MR}{\relax\ifhmode\unskip\space\fi MR }
\providecommand{\MRhref}[2]{%
  \href{http://www.ams.org/mathscinet-getitem?mr=#1}{#2}
}
\providecommand{\href}[2]{#2}
\begin{thebibliography}{10}

\bibitem{AL}
S.~Akiyama and B.~Loridant, \emph{Boundary parametrization of self-affine
  sets}, J. Math. Soc. Japan \textbf{63} (2011), no.~2, 525--579.

\bibitem{BaaFG}
M.~Baake, D.~Frettl\"oh, and U.~Grimm, \emph{A radial analogue of poisson's
  summation formula with applications to powder diffraction and pinwheel
  patterns}, J. Geom. Physics \textbf{57} (2007), 1331--1343.

\bibitem{Baake2013}
M.~Baake and U.~Grimm, \emph{Aperiodic order, vol. 1: A mathematical
  invitation}, Cambridge University Press, Cambridge, 2013.

\bibitem{Ba5}
C.~Bandt, \emph{Self-similar sets 5. integer matrices and fractal tilings of
  ${\mathbb r}^n$}, Proc. Amer. Math. Soc. \textbf{112} (1991), 549--562.

\bibitem{Ba97}
\bysame, \emph{Self-similar tilings and patterns described by mappings}, The
  Mathematics of Long-range Aperiodic Order (R.V. Moody, ed.), NATO ASI Series,
  vol. C 489, Kluwer Academic Publishers, 1997, pp.~45--84.

\bibitem{BG}
C.~Bandt and S.~Graf, \emph{Self-similar sets 7. a characterization of
  self-similar fractals with positive hausdorff measure}, Proc. Amer. Math.
  Soc. \textbf{114} (1992), 995--1001.

\bibitem{BM}
C.~Bandt and M.~Mesing, \emph{Self-affine fractals of finite type}, Convex and
  fractal geometry, Banach Center Publ., vol.~84, Polish Acad. Sci. Inst.
  Math., Warsaw, 2009, pp.~131--148.

\bibitem{Bar}
M.~F. Barnsley, \emph{Fractals everywhere}, 2 ed., Academic Press, 1993.

\bibitem{CR}
J.H. Conway and C.~Radin, \emph{Quaquaversal tilings and rotations}, Invent.
  Math. \textbf{132} (1998), 179--188.

\bibitem{DKV}
P.~Duvall, J.~Keesling, and A.~Vince, \emph{The hausdorff dimension of the
  boundary of a self-similar tile}, J. London Math. Soc. \textbf{61} (2000),
  748--760.

\bibitem{Fal}
K.~J. Falconer, \emph{Fractal geometry: mathematical foundations and
  applications}, 3 ed., J. Wiley \& sons, 2014.

\bibitem{FW}
N.P. Frank and M.F. Whittaker, \emph{A fractal version of the pinwheel tiling},
  Math. Intellig. \textbf{33} (2011), 7--17.

\bibitem{F}
D.~Frettl\"{o}h, \emph{Substitution tilings with statistical circular
  symmetry}, European J. Comb. \textbf{29} (2008), 1881--1893.

\bibitem{GG}
G.~Gelbrich, \emph{Crystallographic reptiles}, Geometria Dedicata \textbf{51}
  (1994), 235--256.

\bibitem{Gi}
W.J. Gilbert, \emph{The fractal dimension of sets derived from complex bases},
  Canad. Math. Bull. \textbf{29} (1986), no.~4, 495--500.

\bibitem{Go}
C.~Goodman-Strauss, \emph{Matching rules and substitution tilings}, Annals
  Math. \textbf{147} (1998), 181--223.

\bibitem{GS}
B.~Gr\"{u}nbaum and G.C. Shephard, \emph{Patterns and tilings}, Freeman, New
  York, 1987.

\bibitem{HLR}
X.-G. He, K.-S. Lau, and H.~Rao, \emph{Self affine sets and graph-directed
  systems}, Constr. Approx. \textbf{19} (2003), 373--397.

\bibitem{LW}
J.C. Lagarias and Y.~Wang, \emph{Integral self-affine tiles in ${\mathbb r}^n.$
  i. standard and non-standard digit sets}, J. London Math. Soc. \textbf{54}
  (1996), 161--179.

\bibitem{LLR}
C.K. Lai, K.-S. Lau, and H.~Rao, \emph{Classification of tile digit sets as
  product-forms}, arXiv 1305.0202.

\bibitem{L}
B.~Loridant, \emph{Crystallographic number systems}, Monatsh. Math.
  \textbf{167} (2012), 511--529.

\bibitem{M}
D.~Mekhontsev, \emph{Ifs tile finder}, \url{https://ifstile.com}.

\bibitem{MPS}
R.V. Moody, D.~Postnikoff, and N.~Strungaru, \emph{Circular symmetry of
  pinwheel diffraction}, Annales Henri Poincar\'e \textbf{7} (2006), 711--730.

\bibitem{R0}
C.~Radin, \emph{Miles of tiles}, Student Mathematical Library, Amer. Math.
  Soc., Providence.

\bibitem{R}
\bysame, \emph{The pinwheel tilings of the plane}, Annals Math. \textbf{139}
  (1994), 661--702.

\bibitem{SaB}
L.~Sadun, \emph{Topology of tiling spaces}, University Lecture Series, Amer.
  Math. Soc., Providence.

\bibitem{Sa}
\bysame, \emph{Some generalizations of the pinwheel tiling}, Discrete Comput.
  Geom. \textbf{20} (1998), 79--110.

\bibitem{ScT}
K.~Scheicher and J.M. Thuswaldner, \emph{Neighbors of self-affine tiles in
  lattice tilings}, Fractals in Graz 2001 (P.~Grabner and W.~Woess, eds.),
  Birkh\"auser, 2003, pp.~241--262.

\bibitem{Se}
M.~Senechal, \emph{Quasicrystals and geometry}, Cambridge University Press,
  Cambridge, 1995.

\bibitem{SoD}
B.~Solomyak, \emph{Dynamics of self-similar tilings}, Ergodic Theory Dyn
  Systems \textbf{17} (1997), 695--738.

\bibitem{So}
\bysame, \emph{Nonperiodicity implies unique composition for self-similar
  translationally finite tilings}, Discrete Comput. Geom. \textbf{20} (1998),
  265--279.

\bibitem{SW}
R.S. Strichartz and Y.~Wang, \emph{Geometry of self-affine tiles 1}, Indiana
  Univ. Math. J. \textbf{48} (1999), 1--24.

\bibitem{Ventrella2012}
J.~Ventrella, \emph{Brainfilling curves - a fractal bestiary -}, Lulu.com,
  Raleigh, North Carolina, 2012, see www.fractalcurves.com.

\end{thebibliography}
\vspace{3ex}

\noindent
Christoph Bandt\\  Institute of Mathematics, University of Greifswald, 17487 Greifswald, Germany. \\ 
\url{bandt@uni-greifswald.de}\vspace{1ex}\\
Dmitry Mekhontsev\\Sobolev Institute of Mathematics, 
4 Acad. Koptyug avenue, 630090 Novosibirsk Russia\\
\url{mekhontsev@gmail.com}\vspace{1ex}\\
Andrei Tetenov\\ Gorno-Altaisk\\
\url{a.tetenov@gmail.com}
\end{document}